\documentclass{amsart}

\usepackage{epsfig}
\usepackage{graphicx}
\usepackage{color}
\usepackage{amsthm,amsmath,amsfonts,amssymb,amscd}
\usepackage{epstopdf}

\usepackage{enumerate}

\def\KK{{\sf {K}}}

\def\1{{\sf {1}}}
\def\0{{\sf {0}}}

\def\n{{\sf {p}}}

\def\qed{\hfill {\hbox{\footnotesize{$\Box$}}}}
\def\0{{\sf {\sf 0}}}

\def\0{{\sf 0}}

\def\so{{\sf supp}}

\def\proof{{\noindent\bf Proof.}\hskip 0.3truecm}
\def\BBox{\kern  -0.2cm\hbox{\vrule width 0.15cm height 0.3cm}}

\def\B{\mathcal{B}}
\def\E{\mathcal{E}}
\def\C{\mathcal{C}}

\def\F{\mathcal{F}}
\def\G{\mathcal{G}}

\def\K{\mathcal{K}}
\def\L{\mathcal{L}}
\def\M{\mathcal{M}}
\def\P{\mathcal{P}}

\def\NN{\mathbb{N}}
\def\ZZ{\mathbb{Z}}

\def\RR{\mathbb{R}}
\def\CC{\mathbb{C}}
\def\KK{\mathbb{K}}

\def\func#1#2{ \colon {#1} \longrightarrow {#2}}

\def\BBox{\kern  -0.2cm\hbox{\vrule width 0.15cm height 0.3cm}}

\oddsidemargin .cm \evensidemargin 1.cm \textwidth=16.5cm
\include {mak}
\begin{document}

\newtheorem{propo}{Proposition}[section]
\newtheorem{lemma}[propo]{Lemma}
\newtheorem{theorem}[propo]{Theorem}
\newtheorem{corollary}[propo]{Corollary}
\newtheorem{definition}[propo]{Definition}

\def\1{{\sf 1}}
\def\n{{\sf n}}
\def\so{{\sf supp}}
\def\proof{{\noindent  Proof.}\hskip 0.3truecm}
\def\BBox{\kern  -0.2cm\hbox{\vrule width 0.15cm height 0.3cm}}
\def\L{{\mathcal  L}}\def\P{{\mathcal  P}}
\def\B{{\mathcal B}}
\def\C{{\mathcal C}}
\def\K{{\mathcal K}}\def\M{{\mathcal M}}
\def\E{{\mathcal E}}
\def\F{{\mathcal F}}
\def\G{{\mathcal G}}
\oddsidemargin 16.5mm
\evensidemargin 16.5mm

\thispagestyle{plain}

%


%

\vspace{5cc}
\begin{center}

{\Large\bf  FLOQUET THEORY FOR  SECOND ORDER LINEAR\\
   HOMOGENEOUS DIFFERENCE EQUATIONS
\rule{0mm}{6mm}\renewcommand{\thefootnote}{}
\footnotetext{The final version will be published in Journal of Difference Equations and Applications\\[1ex]
\scriptsize 2000 Mathematics Subject Classification: 39A10, 11B, 33C45\\
Keywords and Prases: Difference equations, Floquet theory, periodic sequences, Chebyshev polynomials. }}

\vspace{1cc} {\large\it A.M. Encinas and M.J. Jim\'enez}

\vspace{1cc}
\parbox{24cc}{{\small{\bf Abstract.}
In this paper we provide  a version of the Floquet's theorem to be applied to any second order difference equations with quasi-periodic coefficients. To do this we extend to second order linear difference equations with quasi-periodic coefficients, the known equivalence between the Chebyshev equations and the  second order linear difference equations with constant coefficients.  So, any second order linear difference equations with quasi-periodic coefficients is essentially equivalent to a Chebyshev equation, whose parameter only depends on the values of the quasi-periodic coefficients and can be determined by a non-linear recurrence. Moreover, we solve this recurrence and obtaining 
a closed expression for this parameter. As a by-product we also obtain a Floquet's type result; that is, the necessary and sufficient condition for the equation has quasi-periodic solutions.}}
\end{center}
\vspace{1cc}

%
%

\vspace{1.5cc}

\section{Introduction and Preliminaries}
\vspace{.5cm}
It is well known the important role that linear homogeneous difference equations play in several problems of the engineering or of the science. However, whereas the expression of the solutions when the coefficients of the equations are constant is widely known,  the same does not happen for variable coefficients, except for the simplest case of first order equations. In \cite{M97} a complete closed form solution of a second order linear homogeneous difference equation with variable coefficients was presented. The solutions are then expressed solely in terms of the given coefficients. 

We are here interested in the Floquet Theory for second order linear homogeneous difference equation; that is, for difference equations whose coefficients are periodic sequences. The main question in this framework is to find the conditions under which the given equation has periodic solutions with the same period than the coefficients, see for instance \cite{A00}. The strategy we follow in this work is to extend the equivalence between second order equations with constant coefficients and Chebyshev equations to the case of periodic coefficients. Then, the characterization of the existence of periodic solutions can be reduced to the same question about Chebyshev equations. We remark that whereas  the equivalence between a second order equation with periodic coefficients and a Chebyshev equation can be established, more or less easily, by induction, the determination of the parameter of the equivalent Chebysev equation is more difficult since it involves a highly non-linear recurrence. In this work we solve this recurrence obtaining a nice closed formula for this parameter in terms of the coefficients of the considered second order equation. As an straightforward consequence, we  obtain the necessary and sufficient condition for the existence of periodic solutions.
\vspace{.5cc}

Throughout  the paper, $\KK$ denotes either $\RR$ or $\CC$, $\ZZ$ the set of integers, $\NN$ the set of nonnegative integers, $\KK^*=\KK\setminus \{0\}$ and $\NN^*=\NN\setminus\{0\}$. A sequence of elements of $\KK$ is a function $z\func{\ZZ}{\KK}$ and we denote by $\ell(\KK)$ the space of all sequences  of elements of $\KK$ and by $\ell(\KK^*)\subset \ell(\KK)$ the subset of the sequences $z\in \ell(\KK)$ such that $z(k)\not=0$ for all $k\in \ZZ$. The null sequence is denoted by $\0$.

Given $z\in \ell(\KK)$ and $p\in \NN^*$, for any $m\in \ZZ$ we denote by $z_{p,m}\in \ell(\KK)$ the subsequence of $z$ defined as
$$z_{p,m}(k)=z(kp+m),\hspace{.25cm}k\in \ZZ.$$
Clearly, any sequence $z\in \ell(\KK)$ is completely determined by the values of the sequences $z_{p,j}$, for $0\le j\le p-1$. In particular, $z_{1,0}=z$, whereas $z_{2,0}$ and $z_{2,1}$ are the subsequences of even or odd indexes, respectively. Moreover, the sequences $z_{1,m}$ are the {\it shift subsequences} of $z$, since $z_{1,m}(k)=z(k+m)$ for any $k\in \ZZ$. Notice that if we also allow $p=-1$, then $z_{-1,m}$ are the {\it flipped shift subsequences} of $z$, since $z_{-1,m}(k)=z(m-k)$ for any $k\in \ZZ$.

The sequence $z\in \ell(\KK)$ is called {\it quasi-periodic  with period $p\in \NN^*$ and ratio $r\in \KK^*$} if  it satisfies that 
$$z(p+k)=r\, z(k), \hspace{.25cm}k\in \ZZ,$$
which also implies that $z(kp+m)=r^kz(m)$ for any $k,m\in \ZZ$. 

Clearly a sequence $z\in \ell(\KK)$ is periodic with period $p$ iff it is quasi-periodic with period $p$ and ratio $r=1$. The set of quasi-periodic sequences with period $p$ and ratio $r$ is denoted by $\ell(\KK;p,r)$ and we define $\ell(\KK^*;p,r)=\ell(\KK;p,r)\cap \ell(\KK^*)$.
Then, $\ell(\KK;p,1)$ consists of the periodic sequences with period  $p$, whereas $\ell(\KK;1,r)$ consists of the geometric sequences with common ratio $r$; that is, if $z\in \ell(\KK;1,r)$, then  $z(k)=z(0)r^k$. In particular $\ell(\KK;1,1)$ consists of all constant sequences and it is identified with $\KK$.

In the sequel we omit the parameter $r$ when it equals $1$.
Therefore, the space of periodic sequences with period $p$ is
denoted simply by $\ell(\KK;p)$ and hence, $\ell(\KK;1)$ consists of
the constant sequences.

 If $z\in \ell(\KK;p,r)$ is not the null sequence, then $r=z(k_0)^{-1}z(k_0+p)$, where  $k_0=\min\{k\in \NN: z(k)\not=0\}$. Therefore, if $z$ is a non-null
quasi-periodic sequence of period $p$, then $z$  is determined by
the $p+1$ values $z(j)$, $j=0,\ldots,p-1$ and $r$ or equivalently
by the values $z(j)$, $j=0,\ldots,p$.

\begin{lemma}
\label{qp:cha}
Given $p\in \NN^*$ and $r\in \KK^*$, then $z\in \ell(\KK;p,r)$ iff $z_{p,m}\in \ell(\KK;1,r)$ for any $m\in \ZZ$. Moreover, $\ell(\KK;p,r)\subset \ell(\KK;np,r^n)$ for any $n\in \NN^*$.\end{lemma}

Given three sequences $a,c\in \ell(\KK^*)$ and $b\in \ell(\KK)$, we
can consider the irreducible homogeneous linear second order difference equation
\begin{equation}
\label{equation}
a(k)z(k+1)-b(k)z(k)+c(k-1)z(k-1)=0,\hspace{.5cm}k\in \ZZ.
\end{equation}
The sequences  $a,b$ and $c$ are called the {\it coefficients of the
Equation \eqref{equation}} and any sequence $z\in \ell(\KK)$
satisfying the Identity \eqref{equation} is called a {\it solution of
the equation}.  It is well-known that for any $z_0,z_1\in \KK$ and any $m\in \ZZ$, 
there exists a unique solution of Equation \eqref{equation}
satisfying $z(m)=z_0$ and $z(m+1)=z_1$. In addition, when $a,b,c\in \ell(\RR)$, then a solution of Equation \eqref{equation} satisfies that $z\in \ell(\RR)$ iff $z(m),z(m+1)\in \RR$ for some $m\in \ZZ$.

The Equation \eqref{equation} has {\it constant coefficients} when $a,c\in \KK^*$ and $b\in \KK$. Linear difference equations with constant coefficients can be characterized as those  satisfying that $z\in \ell(\KK)$ is a solution iff any shift of $z$ is also a solution. Moreover, $a=c$ iff, in addition, when $z\in \ell(\KK)$ is a solution any flipped shift of $z$ is also a solution. We conclude this section with a result about quasi periodic solutions of equations with constant coefficients.

\begin{lemma}
\label{sol:per}
Given $a,c\in \KK^*$, $b\in \KK$, and $z\in \ell(\KK)$ a solution of the difference equation 
$$az(k+1)-bz(k)+cz(k-1)=0,\hspace{.5cm}k\in \ZZ,$$
then $z\in \ell(\KK;p,r)$ iff $z(p)=r\,z(0)$ and $z(p+1)=r\,z(1)$.
\end{lemma}

\section{Chebyshev sequences}

Among equations with constant coefficients, the so-called {\it Chebyshev equations}  play a main role in the study of this kind of difference equations, see \cite[Theorem 3.1]{ABD05}.
Given $q\in \KK$, the second order difference equation with constant coefficients
\begin{equation}
\label{Chebyshev:equation}
z(k+1)-2qz(k)+z(k-1)=0,\hspace{.25cm} k\in \ZZ,
\end{equation}
is called {\it Chebyshev equation with parameter $q$} and its  solutions are called {\it Chebyshev sequences with parameter $q$}. Clearly any shift of any flipped shift of a Chebyshev sequence is also a Chebyshev sequence. Moreover, a non null Chebyshev sequence determines its parameter, since if $z\in \ell(\KK)$ is a Chebyshev sequence with parameters $q$ and $\hat q$, then 
$$2qz(k)=z(k+1)+z(k-1)=2\hat qz(k), \hspace{.25cm}k\in \ZZ,$$
and hence, $2(q-\hat q)z=\0$, which implies that $q=\hat q$.

Recall that a polynomial sequence $\{P_k(x)\}_{k\in \ZZ}\subset \CC[x]$ is   a {\it sequence of  Chebyshev polynomials} if it satisfies the following three-term recurrence
\begin{equation}
\label{Ch:poly}
P_{k+1}(x)=2xP_k(x)-P_{k-1}(x),\hspace{.25cm}k\in \ZZ.
\end{equation}
Therefore, any sequence of  Chebyshev polynomials is completely determined by the choice of the two polynomials $P_{m}(x)$ and $P_{m+1}(x)$ for some $m\in \ZZ$. In particular, the choice $U_{-1}(x)=0$ and $U_0(x)=1$, determines the sequence $\{U_k(x)\}_{k\in \ZZ}$, where for any $k\in \ZZ$, $U_k(x)$   is the {\it $k$-th Chebyshev
polynomial of second kind}, see \cite{MH03}. Then $U_{-k}(x)=-U_{k-2}(x)$, $U_k(-x)=(-1)^kU_k(x)$ for any $k\in \ZZ$ and any $x\in \CC$  and moreover  
\begin{equation}
\label{ch:second}
U_k(x)=\sum\limits_{j=0}^{\lfloor \frac{k}{2}\rfloor}(-1)^j{k-j\choose j}(2x)^{k-2j},\hspace{.25cm}k\in \NN.
\end{equation}

In addition, for any sequence of Chebyshev polynomials $\{P_k(x)\}_{k\in \ZZ}$ we have that 
\begin{equation}
\label{comb}
P_k(x)=P_0(x)U_k(x)-P_{-1}(x)U_{k-1}(x),\hspace{.25cm}\hbox{for any}\hspace{.25cm}k\in \ZZ. 
\end{equation}
In particular,  $\{P_k(x)\}_{k\in \ZZ}\subset \RR[x]$ iff $P_{-1}(x),P_0(x)\in \RR[x]$,  since $U_k(x)\in \RR[x]$ for any $k\in \ZZ$.   Moreover, if given $p\in \NN^*$ we apply \eqref{comb} to the flipped shift sequence $\{U_{p-k}(x)\}_{k\in \ZZ}$ we obtain the well-known identity
\begin{equation}
\label{fs}
U_{p-k}(x)=U_{p}(x)U_k(x)-U_{p+1}(x)U_{k-1}(x),\hspace{.25cm}\hbox{for any}\hspace{.25cm}k\in \ZZ. 
\end{equation}

The sequences $\{T_k(x)\}_{k\in \ZZ}$, $\{V_k(x)\}_{k\in \ZZ}$ and $\{W_k(x)\}_{k\in \ZZ}$ defined by taking $T_{-1}(x)=x$, $V_{-1}(x)=1$, $W_{-1}(x)=-1$ and $T_0(x)=V_0(x)=W_0(x)=1$ are known as {\it Chebyshev polynomials of first, third and fourth order}, respectively.
From the Identity \eqref{comb} we obtain that $V_k(x)=U_k(x)-U_{k-1}(x)$, $W_k(x)=U_k(x)+U_{k-1}(x)$ $T_k(x)=U_k(x)-xU_{k-1}(x)=\frac{1}{2}\big[U_k(x)-U_{k-2}(x)\big]$ for any $k\in \ZZ$, which   implies that  $T_{-k}(x)=T_{k}(x)$, $T_k(-x)=(-1)^kT_k(x)$ for any $k\in \ZZ$ and any $x\in \CC$ and moreover 
\begin{equation}
\label{ch:first}
T_k(x)=\dfrac{k}{2}\sum\limits_{j=0}^{\lfloor \frac{k}{2}\rfloor}\dfrac{(-1)^j}{k-j}{k-j\choose j}(2x)^{k-2j},\hspace{.25cm}k\in \NN^*.
\end{equation}
 \vspace{.5cc}

Backing to Chebyshev sequences, it is clear that any Chebyshev sequence with parameter $q$ is of the form $\{P_k(q)\}_{k\in \ZZ}$, where $\{P_k(x)\}_{k\in \ZZ}$ is a sequence of Chebyshev polynomials. Therefore, many properties of Chebyshev sequences are consequence of properties of Chebyshev polynomials and conversely.

\begin{propo}
\label{Cheb:periodic} Given $q\in \KK$, then for any $p\in \NN^*$  the Chebyshev equation with parameter $q$   has non-null solutions belonging to $\ell(\KK;p,r)$  iff 
$$r=T_p(q)\pm \sqrt{T_p(q)^2-1}.$$
In particular,   the following results hold:
\begin{itemize}
\item[(i)] If $p\in \NN^*$, the Chebyshev equation with parameter $q$ has non-null solutions belonging to $\ell(\KK;p)$ iff 
${q=\cos\big(\frac{2j\pi}{p}\big)}$, $j=0,\ldots,\lceil\frac{p-1}{2}\rceil$. 
%
\item[(ii)] If $r\in \KK^*$ the Chebyshev equation with parameter $q$ has non-null solutions belonging to $\ell(\KK;1,r)$ iff 
$q=\dfrac{1}{2}(r+r^{-1})$.
\item[(iii)] The Chebyshev equation with parameter $q$ has constant solutions iff $q=1$.
\end{itemize} 
\end{propo}
\proof Given  $z(k)=AU_k(q)+BU_{k-1}(q)$ a non-null Chebyshev sequence with parameter $q$,  then  from Lemma \ref{sol:per} $z\in \ell(\KK;p,r)$ iff $z(p)=rz(0)$ and $z(p+1)=rz(1)$; that is iff 
$$\begin{bmatrix}U_p(q)-r &U_{p-1}(q)\\[1ex]
U_{p+1}(q)-2qr & U_p(q)-r\end{bmatrix}\begin{bmatrix}A\\[1ex]B\end{bmatrix}=\begin{bmatrix}0\\[1ex]0\end{bmatrix}$$
Therefore, the Chebyshev equation with parameter $q$ has non-null solutions belonging to $\ell(\KK;p,r)$ iff the determinant of the above matrix equals $0$; that is, applying \eqref{fs}, iff 
$$0= r^2+U_p^2(q)-U_{p+1}(q)U_{p-1}(q)-2r\big[U_{p}(q)-qU_{p-1}(q)\big]+1=r^2-2rT_p(q)+1.$$

(i) When $z\in \ell(\KK;p)$; that is, when $r=1$, the above equation becomes $T_p(q)=1$, that implies that $q=\cos\Big(\frac{2j\pi}{p}\Big)$, $j=0,\ldots,\lceil\frac{p-1}{2}\rceil$, see  \cite{MH03}. 

 (ii) When $p=1$, then the above equation becomes $r^2-2qr+1=0$ and hence $q=\dfrac{1}{2}(r+r^{-1})$.

(iii)  It is an straightforward consequence of both, (i) or (ii). \qed
\vspace{.5cc}

The role of Chebyshev equations among constant coefficients equations is described by the following results. We start with an easy-to-proof  result involving first order linear difference equations with constant coefficients.
\begin{lemma} 
\label{order:1} Let $r\in \KK^*$ and consider
$q=\frac{1}{2}(r+r^{-1})$. Then, a sequence $z\in
\ell(\KK)$ is a solution of the first order difference equation
$z(k+1)=r  z(k)$, $k\in \ZZ$, iff it is a  solution of the
Chebyshev equation  $z(k+1)-2q z(k)+z(k-1)=0$, $k\in \ZZ$, satisfying
that $z(1)=r z (0)$ or equivalently, iff $z$ is a multiple of $\alpha U_{k-1}(q)-U_{k-2}(q)$.\end{lemma}

The second result concerns to the even and odd subsequences of a given Chebyshev sequence.

\begin{lemma}
\label{Ch:odd-even}
Let  $q\in \KK$ and $z\in \ell(\KK)$ be a solution of the Chebyshev equation
$$z(k+1)-2qz(k)+z(k-1)=0,\hspace{.25cm} k\in \ZZ.$$
Then for any $m\in \ZZ$, the subsequence  $z_{2,m}$ is a  Chebyshev sequence  with parameter $2q^2-1$.
\end{lemma}
\proof When $q=0$, then $z(k)=-z(k-2)$ for any $k\in \ZZ$; which implies that $z_{2,m}(k)=-z_{2,m}(k-1)$, for any $k\in \ZZ$. Applying Lemma \ref{order:1} we obtain that $z_{2,m}$ is solution of the Chebyshev equation with parameter $-1=2q^2-1$. 

When $q\in \KK^*$,  for any $k\in \ZZ$ we have that  $z(k)=\dfrac{1}{2q}\big[z(k+1)+z(k-1)\big]$  and hence, 
$$\begin{array}{rl}
0=&\hspace{-.25cm}z(2k+m+1)-2qz(2k+m)+z(2k+m-1)\\[1ex]
=&\hspace{-.25cm}\dfrac{1}{2q}\big[z(2k+m+2)+2z(2k+m)+z(2k+m-2)\big]-2qz(2k+m)\\[2ex]
=&\hspace{-.25cm}\dfrac{1}{2q}\big[z(2k+m+2)-2(2q^2-1)z(2k+m)+z(2k+m-2)\big]\\[2ex]
=&\hspace{-.25cm}\dfrac{1}{2q}\big[z_{2,m}(k+1)-2(2q^2-1)z_{2,m}(k)+z_{2,m}(k-1)\big]\end{array}$$
and the  claim follows.\qed

Many properties and identities involving  Chebyshev polynomials are consequence of being solutions of the Chebyshev equations. For instance, from the above lemma  we have the following classical identities, see $(1.14)$ and $(1.15)$ in \cite{MH03},
\begin{equation}
\label{doubling:second}
U_{2k}(x)=W_k(2x^2-1)\hspace{.25cm}\hbox{and}\hspace{.25cm}U_{2k+1}(x)=2x U_k(2x^2-1), \hspace{.25cm}k\in \ZZ,
\end{equation}
which in turns implies that 
\begin{equation}
\label{doubling:first}
T_{2k}(x)=T_k(2x^2-1)\hspace{.25cm}\hbox{and}\hspace{.25cm}T_{2k+1}(x)=xV_k(2x^2-1), \hspace{.25cm}k\in \ZZ.
\end{equation}

The next result shows that any difference equation with constant coefficients is equivalent to a Chebyshev equation. Although it is a known result, see \cite[Theorem 3.1]{ABD05},   we reproduce here its proof, for the sake of completeness. 

\begin{theorem}
\label{Ch:complex}
Consider $a,c\in \KK^*$, $b\in \KK$ and $z\in \ell(\KK)$ a solution of the second order difference equation with constant coefficients
$$az(k+1)-bz(k)+cz(k-1)=0,\hspace{.25cm} k\in \ZZ.$$
Then $z(k)=(\sqrt{a^{-1}c})^{k}v(k)$, $k\in \ZZ$, where $v\in \ell(\CC)$ is a solution of the Chebyshev equation 
$$v(k+1)-2qv(k)+v(k-1)=0,\hspace{.25cm} k\in \ZZ$$
whose parameter is $q=\dfrac{b}{2\sqrt{ac}}$. Moreover, if $\KK=\RR$ and $ac>0$ then $v \in\ell(\RR)$, whereas when $ac<0$ then, for any $m\in \ZZ$, $z_{2,m}(k)=( a^{-1}c )^{k}w(k)$, where $w\in \ell(\RR)$ is a solution of the Chebyshev equation 
$$w(k+2)-2\hat qw(k+1)+w(k)=0,\hspace{.25cm} k\in \ZZ$$
whose parameter is $\hat q=\dfrac{b^2}{2ac}-1\in \RR$. 
\end{theorem}
\proof Clearly, for any $k\in \ZZ$ we have that 
$$0=az(k+1)-bz(k)+cz(k-1)=(\sqrt{a^{-1}c})^{k-1}\Big[cv(k+1)-b\sqrt{a^{-1}c}\,v(k)+cv(k-1)\Big].$$
Therefore, $v\in \ell(\CC)$ is a solution of the difference equation with complex coefficients 
$$ 0=cv(k+1)-b\sqrt{a^{-1}c}v(k)+cv(k-1)$$
or, equivalently of the Chebyshev equation with parameter $q=\dfrac{b\sqrt{a^{-1}c}}{2c}=
\dfrac{b}{2\sqrt{ac}}$. When $\KK=\RR$,  if $ac>0$, then $q\in \RR$ and moreover $v\in \ell(\RR)$; whereas if $ac<0$, then  $\hat q=2q^2-1=\dfrac{b^2}{2ac}-1\in \RR$. Moreover, given $m\in \ZZ$, $(\sqrt{a^{-1}c})^mv$ is also a solution of the Chebyshev equation with parameter $q$ and hence applying Lemma \ref{Ch:odd-even}, $w=(\sqrt{a^{-1}c})^mv_{2,m}$ is a solution of the Chebyshev equation with parameter $\hat q\in \RR$. Therefore,   for any $k\in \ZZ$ we have
$$ w(k)= (\sqrt{a^{-1}c})^mv(2k+m)=(\sqrt{a^{-1}c})^m (\sqrt{a^{-1}c})^{-m-2k}z(2k+m)=(a^{-1}c)^{-k}z_{2,m}(k)$$
which in particular implies that $w\in \ell(\RR)$ and hence the result. \qed
\vspace{.5cc}

No we can derive a Floquet's type theorem for equations with constant coefficients.
\begin{corollary}
\label{floquet:constant}
Given $a,c\in \KK^*$ and $b\in \KK$, the equation with constant coefficients 
 $$az(k+1)-bz(k)+cz(k-1)=0,\hspace{.5cm}k\in \ZZ$$
has  quasi-periodic solutions of period $p\in \NN^*$ and ratio $r\in\KK^*$ iff 
$$r=\sqrt{\dfrac{c^p}{a^p}}\left[T_p(q)\pm \sqrt{T_p(q)^2-1}\right],\hspace{.25cm}\hbox{where}\hspace{.25cm}q=\dfrac{b}{2\sqrt{ac}}.$$
Therefore, the equation has geometric solutions with ratio $r$ iff $r=\dfrac{b\pm \sqrt{b^2-4ac}}{2a}$ and, in particular, it has constant solutions iff $b=a+c$.
\end{corollary}
\proof According to Theorem \ref{Ch:complex}, the $z\in \ell(\KK)$ is a solution of the above equation iff the sequence  defined for any $k\in \ZZ$ as $v(k)=\alpha^kz(k)$, where $\alpha=\sqrt{ac^{-1}}$, is a  Chebyshev sequence with parameter $q=\dfrac{b}{2\sqrt{ac}}$. Therefore, $z\in \ell(\KK;p,r)$ iff $v\in \ell(\KK;p,r\alpha^p)$ and from Proposition \ref{Cheb:periodic}, this happens iff 
$$ r\alpha^p=T_p(q)\pm \sqrt{T_p(q)^2-1}.$$
\qed
   \vspace{.5cc}

Our aim in this paper is to extend the above results to a wider class of linear difference equations.

\section{Second Order Difference Equations with Quasi-Periodic Coefficients}

We say that the Equation \eqref{equation} has {\it quasi-periodic
coefficients with period $p\in \NN^*$ and ratio $r\in \KK^*$}
if $a,c\in \ell(\KK^*;p,r)$ and $b\in\ell(\KK;p,r)$. In
particular, we say that the Equation \eqref{equation} has {\it
constant coefficients} when $a,c\in \ell(\KK^*;1,1)\equiv \KK^*$ and $b\in \ell(\KK;1,1)\equiv \KK$.  

The Equation \eqref{equation} is called {\it symmetric} when $a=c$. When $a,b,c\in \ell(\RR)$, symmetric equations are also called {\it self-adjoint equations}.  It
is well-known that any irreducible second order linear difference
equation is equivalent to a symmetric one, and to a self-adjoint one when $a,b,c\in \ell(\RR)$. With  this end, we consider  
%
the function $\phi\func{\ell(\KK^*)\times \ell(\KK^*)}{\ell(\KK^*)}$ defined as
\begin{equation}
\label{sym}\phi(a,c)(0)=1,\hspace{.25cm}\phi(a,c)(k)=\displaystyle
\prod\limits_{j=0}^{k-1}\dfrac{a(j)}{c(j)},\hspace{.25cm}\hbox{when $k>0$}\hspace{.25cm}\hbox{and}\hspace{.25cm}\phi(a,c)(k)=\displaystyle
\prod\limits_{j=k}^{-1}\dfrac{c(j)}{a(j)},\hspace{.25cm}\hbox{when $k<0$}.
\end{equation}
 
\begin{lemma}
\label{equivalence} Given $a,c\in \ell(\KK^*)$, the following
properties hold:
\begin{itemize}
\item[(i)] $\phi(a,c)(k)=1$ for all $k\in \ZZ$ iff $a=c$.
\item [(ii)] $\phi(a,c)(k-1)a(k-1)=\phi(a,c)(k)c(k-1)$, $k\in \ZZ.$ Therefore,
$z\in \ell(\KK)$ is a solution of the difference equation whose
coefficients are $a,b$ and $c$ iff it is a solution of the symmetric  difference equation whose coefficients are $\phi(a,c) a$ and
$\phi(a,c) b$.
\item[(iii)] If $a,c\in \ell(\KK^*;p,r)$, then
$\phi(a,c)\in \ell\big(\KK^*;p,\phi(a;c)(p)\big)$. Therefore, if $b\in \ell(\KK;p,r)$, then 
$\phi(a,c)b\in \ell\big(\KK^*;p,r\phi(a;c)(p)\big)$.
\end{itemize}
\end{lemma}
%

The next result shows the role that Chebyshev equations play to solve some linear systems of difference equations with constant
coefficients.  It is the key to solve general second order linear
difference equations with quasi-periodic coefficients.

\begin{propo}
\label{complex:system} Given $p\in \NN^*$, $a_j\in \KK^*$ and  $b_j\in \KK$, $j=0,\ldots,p-1$,  consider the sequences $v_j\in \ell(\KK)$, $j=0,\ldots,p-1$, satisfying  the equalities
$$\left\{\begin{array}{rll}
  b_0v_0(k)=&\hspace{-.25cm} a_0v_1(k)+ a_{p-1}v_{p-1}(k-1),& \\[1ex]
 b_jv_j(k)=&\hspace{-.25cm} a_jv_{j+1}(k)+ a_{j-1}v_{j-1}(k), & j=1,\ldots,p-2,\\[1ex]
 b_{p-1}v_{p-1}(k)=&\hspace{-.25cm} a_{p-1}v_0(k+1)+ a_{p-2}v_{p-2}(k),&
\end{array}\right.$$
where $v_1(k)=v_0(k+1)$, $k\in \ZZ$,  when $p=1$. Then,  there exists   $q_p(a_0,\ldots,a_{p-1};b_0,\ldots,b_{p-1})\in \KK$ such that for any $j=0,\ldots,p-1$, $v_j$ is a solution of the Chebyshev equation
$$ z(k+1)-2q_p(a_0,\ldots,a_{p-1};b_0,\ldots,b_{p-1})z(k)+z(k-1)=0,\hspace{.25cm}k\in \ZZ.$$
Moreover, if $a_j,b_j\in \RR$, $j=0,\ldots,p-1$, then  for any $j=1,\ldots,p-1$ it is verified that 
$$iq_p(a_0,\ldots,\pm i a_j,\ldots,a_{p-1};b_0,\ldots,b_{p-1})\in \RR.$$
\end{propo}
\proof We prove the claim by induction on $p$.

If $p=1$,  the system is reduced to the equation
$b_0v_0(k)=a_0v_0(k+1)+ a_0v_0(k-1)$
and hence it suffices to take $q_1(a_0;b_0)=\dfrac{b_0}{2a_0}\in \KK$. Moreover, if $a_0,b_0\in \RR$, then   $iq_1(\pm ia_0;b_0)=\pm q_1(a_0;b_0)\in \RR$.

If $p=2$, then the system becomes
$$\left\{\begin{array}{rl}
b_0v_0(k)=&\hspace{-.25cm}a_0v_1(k)+ a_1v_1(k-1),  \\[1ex]
b_1v_1(k)=&\hspace{-.25cm}a_1v_0(k+1)+a_0v_0(k).
\end{array}\right.$$

If $b_1\not=0$,  obtaining  $v_1$ from the second equation and substituting its value at the first one, we get
$$b_0v_0(k)=\dfrac{1}{b_1}\Big[a_0a_1v_0(k+1)+(a_0^2+ a_1^2)v_0(k)+ a_0a_1v_0(k-1)\Big],$$
which implies that $q_2(a_0,a_1;b_0,b_1)=\dfrac{1}{2a_0a_1}\big[b_0b_1-a_0^2- a_1^2\big]$. Moreover, when $a_0,a_1,b_0,b_1\in \RR$, have that 
$$iq_2(\pm ia_0,a_1;b_0,b_1)=\dfrac{\pm1}{2a_0a_1}\big[b_0b_1+a_0^2-a_1^2\big]\in \RR\hspace{.25cm}\hbox{and}\hspace{.25cm}iq_2(a_0,\pm ia_1;b_0,b_1)= \dfrac{\pm 1}{2a_0a_1}\big[b_0b_1-a_0^2+a_1^2\big]\in \RR.$$

As the above Chebyshev equation with parameter $q_2$  has constant coefficients, the sequence $v\in \ell(\CC)$ defined for any $k\in \ZZ$ as $v(k)=v_0(k+1)$  is also a solution of the same equation. Therefore, $v_1$ is also a solution since, from the second equation of the system, it is linear combination of the sequences $v$ and $v_0$.

If  $b_0\not=0$, then obtaining $v_0$ from the first equation and substituting its value at the second one, we get
$$b_1v_1(k)=\dfrac{1}{b_0}\Big[a_0a_1v_1(k+1)+(a_0^2+ a_1^2)v_1(k)+ a_0a_1v_1(k-1)\Big],$$
which newly implies the same conclusions than above.

If $b_0=b_1=0$, then  $v_1(k)= r v_1(k-1)$ and
$v_0(k+1)=r^{-1}v_0(k)$, where $r=-\dfrac{a_1}{a_0}$.
Therefore, applying Lemma \ref{order:1}, we get the both, $v_0$ and $v_1$ are solutions of the Chebyshev equation with parameter 
$\frac{1}{2}(r+r^{-1})=q_2(a_0,a_1;0,0)$. \vspace{.5cc}


  Suppose now that $p\ge 3$ and that the claims are true for any $1\le\ell \le p-1$.

If $b_{p-1}\not=0$, then  from the last equation we have
$$v_{p-1}(k)=b_{p-1}^{-1}a_{p-1}v_0(k+1)+b_{p-1}^{-1}a_{p-2}v_{p-2}(k),\hspace{.5cm}\hbox{for any $k\in \ZZ$}$$
and substituting the value of $v_{p-1}(k-1)$ and of $v_{p-1}(k)$  at the first and at the penultimate equations of the system, we get
$$\left\{\begin{array}{rll}
b_{p-1}^{-1}(b_0b_{p-1}- a_{p-1}^2)v_0(k)=&\hspace{-.25cm}a_0v_1(k)+ b_{p-1}^{-1}a_{p-1}a_{p-2}v_{p-2}(k-1), & \\[1ex]
b_jv_j(k)=&\hspace{-.25cm}a_jv_{j+1}(k)+a_{j-1}v_{j-1}(k),& j=1,\ldots,p-3,\\[1ex]
b_{p-1}^{-1}(b_{p-2}b_{p-1}-a_{p-2}^2)v_{p-2}(k)=&\hspace{-.25cm}b_{p-1}^{-1}a_{p-2}a_{p-1}v_0(k+1)+a_{p-3}v_{p-3}(k).&
\end{array}\right.$$
Applying the induction hypothesis and taken
$$ q_p= 
q_{p-1}\left(a_0,a_1,\ldots,a_{p-3},b_{p-1}^{-1}a_{p-2}a_{p-1};b_{p-1}^{-1}(b_0b_{p-1}-a_{p-1}^2),b_1,\ldots,b_{p-3},b_{p-1}^{-1}(b_{p-2}b_{p-1}-a_{p-2}^2)\right)\in \KK,
$$
for any $j=0,\ldots,p-2$ it is satisfied that
$$2q_pv_j(k)= v_j(k+1)+ v_j(k-1),\hspace{.25cm}\hbox{for any $k\in \ZZ$}.$$
Moreover, since $v_{p-1}$ is a linear combination of two solutions of the same Chebyshev equation, it is also a solution of it. Furthermore, applying  the induction hypothesis,   when $a_j,b_j\in \RR$, $j=0,\ldots,p-1$, then $b_{p-1}^{-1}a_{p-2}a_{p-1},b_{p-1}^{-1}(b_0b_{p-1}-a_{p-1}^2),b_{p-1}^{-1}(b_{p-2}b_{p-1}-a_{p-2}^2)\in \RR$ and we can also conclude that 
$$  i q_p(a_0,\ldots,\pm i a_j,\ldots,a_{p-1};b_0,\ldots,b_{p-1})\in \RR,$$ 
for any $j=0,\ldots,p-1$.

When $b_0\not=0$, obtaining $v_0$ from the first equation and applying the same reasoning   than above, for any $j=0,\ldots,p-1$ we get
$$2q_pv_j(k)= v_j(k+1)+v_j(k-1),\hspace{.25cm}\hbox{for any $k\in \ZZ$}.$$
where
$$q_p=q_{p-1}\left(a_1,\ldots,a_{p-2},b_{0}^{-1}a_{0}a_{p-1};b_{0}^{-1}(b_{0}b_{1}-a_{0}^2),b_2,\ldots,b_{p-2},b_{0}^{-1}(b_0b_{p-1}-a_{p-1}^2)\right)\in \KK$$
and the remaining properties for $q_p$ are also true.

If  $b_0=b_{p-1}=0$, then
$v_{0}(k)=-a_{p-1}^{-1}a_{p-2}v_{p-2}(k-1)$ and
$v_{p-1}(k)=-a_{p-1}^{-1}  a_{0} v_{1}(k+1)$, and hence
substituting its values at the second and at the penultimate
equations, we obtain
$$\left\{\begin{array}{rll}
b_1v_1(k)=&\hspace{-.25cm}a_1v_2(k)-a_{p-1}^{-1}a_0a_{p-2} v_{p-2}(k-1), & \\[1ex]
b_jv_j(k)=&\hspace{-.25cm}a_jv_{j+1}(k)+a_{j-1}v_{j-1}(k),  & j=2,\ldots,p-3,\\[1ex]
b_{p-2}v_{p-2}(k)=&\hspace{-.25cm}-a_{p-1}^{-1} 
a_{0}a_{p-2}v_{1}(k+1)+a_{p-3}v_{p-3}(k).&
\end{array}\right.$$

Applying newly the induction hypothesis and taken
$$q_p=q_{p-2}\left(a_1,\ldots,a_{p-3},-a_{p-1}^{-1} a_0a_{p-2};b_1,\ldots,b_{p-2}\right)\in \KK,$$
then for any $j=1,\ldots,p-2$ it is satisfied that
$$2q_pv_j(k)= v_j(k+1)+v_j(k-1),\hspace{.25cm}\hbox{for any $k\in \ZZ$}.$$
Finally, as the sequences $v_0$ and  $v_{p-1}$ are both multiple of the solutions of the above Chebyshev equation, they are also solutions of it.
Moreover,  when $a_j,b_j\in \RR$, $j=0,\ldots,p-1$, it is also clear that  
$$  i q_p(a_0,\ldots,\pm i a_j,\ldots,a_{p-1};b_0,\ldots,b_{p-1})\in \RR,\hspace{.25cm}\hbox{for any $j=1,\ldots,p-1$,}$$
since   then $-a_{p-1}^{-1} a_0a_{p-2}\in \RR$.   \qed
\vspace{.5cc}

\begin{theorem}
\label{Ch:sys} Consider $p\in \NN^*$, $r\in \KK^*$, $a,c\in \ell(\KK^*;p,r)$, $b\in \ell(\KK;p,r)$,
$s=\displaystyle \prod\limits_{j=0}^{p-1}\dfrac{a(j)}{c(j)}$, $\gamma=(\sqrt{rs})^{-1}$ and
$z\in \ell(\KK)$ a solution of the  equation with quasi-periodic
coefficients
$$a(k)z(k+1)-b(k)z(k)+c(k-1)z(k-1)=0,\hspace{.25cm}k\in \ZZ.$$
Then, there exists $q_{p,r}(a;b;c)  \in \CC$ such that  for any $m\in \ZZ$, $z_{p,m}(k)=\gamma^{k}v(k)$, $k\in \ZZ$, where  the sequence $v\in \ell(\CC)$ is a solution of the Chebyshev equation
$$ v(k+1)-2 q_{p,r}(a;b;c)v(k)+v(k-1)=0,\hspace{.25cm}k\in \ZZ.$$
Moreover, if $\KK=\RR$ we have the following results:
\begin{itemize}
\item[(i)] If $rs>0$, then $q_{p,r}(a;b;c)\in \RR$ and $v\in \ell(\RR)$. 
\item[(ii)] If $rs<0$, then $q_{p,r}(a;b;c)^2\in \RR$ and then $z_{2p,m}(k)=(rs)^{-k}u(k)$, $k\in \ZZ$,  where $u\in \ell(\RR)$ is  a solution of the Chebyshev equation
$$ u(k+1)-2 \big(2q_{p,r}(a;b;c)^2-1\big)u(k)+u(k-1)=0,\hspace{.25cm}k\in \ZZ.$$ 
\end{itemize}
\end{theorem}
\proof  Given $m\in \ZZ$, then $m=k_0p+j$, where  $0\le j\le p-1$. Therefore,  $z_{p,m}(k)=z_{p,j}(k+k_0)$  and hence it suffices to prove the claims   for $0\le j\le p-1$ and take into account that if a sequence is a solution of a difference equation with constant coefficients, then any shift  is also a solution of the same equation.

From Part (ii) of Lemma \ref{equivalence}, we know that $z\in \ell(\KK)$ is a solution of the given equation iff it is a solution of the symmetric equation 
$$\phi(a,c)(k)a(k) z(k+1)-\phi(a,c)(k)b(k)z(k)+\phi(a,c)(k-1)a(k-1)z(k-1)=0,\hspace{.25cm}k\in \ZZ$$
and moreover, part (iii) of  Lemma \ref{equivalence} implies that $\phi(a,c)a,\phi(a,c) b\in \ell(\KK;p,rs)$.

Since the coefficients of this last equation are  quasi-periodic with period $p$ and ratio $rs$, $z\in \ell(\KK)$ is a solution iff the subsequences $z_{p,j}$, $j=0,\ldots,p-1$ satisfy the equalities
$$\left\{\begin{array}{rl}
\phi(a,c)(0)b(0)z_{p,0}(k)=&\hspace{-.25cm}\phi(a,c)(0)a(0)z_{p,1}(k)+\gamma^{2}\phi(a,c)(p-1) a(p-1)z_{p,p-1}(k-1),  \\[1ex]
\phi(a,c)(j)b(j)z_{p,j}(k)=&\hspace{-.25cm}\phi(a,c)(j)a(j)z_{p,j+1}(k)+\phi(a,c)(j-1)a(j-1)z_{p,j-1}(k) ,\hspace{.15cm} 1\le j\le p-2,\\[1ex]
\phi(a,c)(p-1)b(p-1)z_{p,p-1}(k)=&\hspace{-.25cm}\phi(a,c)(p-1) a(p-1)z_{p,0}(k+1)+\phi(a,c)(p-2)a(p-2)z_{p,p-2}(k).
\end{array}\right.$$

Defining   $v_j(k)=\gamma^{-k}z_{p,j}(k)$, $j=0,\ldots,p-1$, we get that 
$$\left\{\begin{array}{rl}
\phi(a,c)(0)b(0)v_0(k)=&\hspace{-.25cm}\phi(a,c)(0)a(0)v_1(k)+\gamma\phi(a,c)(p-1) a(p-1)v_{p-1}(k-1)  \\[1ex]
\phi(a,c)(j)b(j)v_j(k)=&\hspace{-.25cm}\phi(a,c)(j)a(j)v_{j+1}(k)+\phi(a,c)(j-1)a(j-1)v_{j-1}(k) ,\hspace{.15cm} 1\le j\le p-2,\\[1ex]
\phi(a,c)(p-1)b(p-1)v_{p-1}(k)=&\hspace{-.25cm}\gamma\phi(a,c)(p-1) a(p-1)v_0(k+1)+ \phi(a,c)(p-2)a(p-2)v_{p-2}(k) 
\end{array}\right.$$
We obtain the result applying Proposition \ref{complex:system} and taking
$$q_{p,r}(a;b;c)=q_{p}\big(\phi(a,c)(0)a(0),\ldots,\gamma\phi(a,c)(p-1)a(p-1); \phi(a,c)(0)b(0),\ldots,\phi(a,c)(p-1)b(p-1)\big).$$

Moreover, it is clear that when $a,b,c\in \ell(\RR)$ and $rs>0$, then $\gamma\in \RR$. Applying newly Proposition \ref{complex:system} we obtain that $q_{p,r}(a;b;c)\in\RR$ and $v\in \ell(\RR)$, which prove (i).

\noindent (ii)  When $rs<0$, then $\gamma=-i(\sqrt{|rs|})^{-1}$. Therefore, if we consider $\hat a\in \ell(\RR)$ defined for any $k,j\in \ZZ$ as $\hat a(pk+j)=a(pk+j)$ if $j\not=p-1$ and as $\hat a(pk+p-1)=(\sqrt{|rs|})^{-1}\,a(pk+p-1)$.
Clearly, 
$$q_{p,r}(a;b;c)=q_p\big(\phi(a,c)(0)\hat a(0),\ldots,-i\phi(a,c)(p-1)\hat a(p-1); \phi(a,c)(0)b(0),\ldots,\phi(a,c)(p-1)b(p-1)\big)$$
that from Proposition \ref{complex:system} implies that $iq_{p,r}(a;b;c)\in \RR$ and hence, $q_{p,r}(a;b;c)^2\in \RR$. The conclusion follows the same reasoning as in the last part of Theorem \ref{Ch:complex}, tacking into account that $z_{p,m}(2k)=z(2pk+m)=z_{2p,m}(k)$, whereas 
$z_{p,m}(2k+1)=z(2pk+p+m)=z_{2p,p+m}(k)$ for any $k\in \ZZ$. 
  \qed

\section{The Floquet Functions}

%

Given $p\in \NN^*$, $r\in \CC$, $a,c\in \ell(\CC^*;p,r)$ and $b\in \ell(\CC,;p,r)$,   Theorem \ref{Ch:sys} establishes that there exists  $q_{p,r}(a;b;c)\in \CC$ such that the difference equation 
$$a(k)z(k+1)-b(k)z(k)+c(k-1)z(k-1)=0,\hspace{.25cm}k\in \ZZ.$$
is equivalent to the Chebyshev equation with parameter $q_{p,r}(a;b;c)$. The aim of this section is to obtain the expression of $q_{p,r}(a;b;c)$. 
Therefore, given $p\in \NN^*$ and $r\in \KK^*$,  we call {\it Floquet function of order $p$ and ratio $r$}, the function $q_{p,r}\func{\ell(\CC^*;p,r)\times \ell(\CC;p,r)\times \ell(\CC^*;p,r)}{\CC}$ such that for any  $a,c\in \ell(\CC^*;p,r)$ and $b\in \ell(\CC;p,r)$,  if $z\in \ell(\CC)$ is a solution of the equation
$$a(k)z(k+1)-b(k)z(k)+c(k-1)z(k-1)=0,\hspace{.25cm}k\in \ZZ,$$
then  for any $m\in \ZZ$, $v(k)=\gamma^{-k}z_{p,m}(k)$ is a solution of the Chebyshev equation
$$ v(k+1)-2q_{p,r}(a;b;c)zv(k)+v(k-1)=0,\hspace{.25cm}k\in \ZZ.$$

Notice  that it suffices to determine the expression of the parameter for the symmetric case and for periodic coefficients. Specifically
if given $a,c\in \ell(\CC^*;p,r)$ and $b\in \ell(\CC;p,r)$ we consider $\gamma=\left(r\prod\limits_{j=0}^{p-1}\dfrac{a(j)}{c(j)}\right)^{-\frac{1}{2}}$ and the pair $(a_\phi,b_\phi)\in \ell(\CC^*;p)\times \ell(\CC;p)$  defined as the periodic extension of 
$$\begin{array}{rlrl}
a_\phi(k)=\phi(a,c)(k)a(k),& k=0,\ldots,p-2; & a_\phi(p-1)=\gamma\phi(a,c)(p-1)a(p-1),\\[1ex]
b_\phi(k)=\phi(a,c)(k)b(k),&  k=0,\ldots,p-1,&\end{array}$$ 
then in the proof of Floquet's Theorem we have shown that $q_{p,r}(a;b;c)=q_{p,1}(a_\phi;b_\phi;a_\phi)$. Moreover, if we consider the function 
$Q_p\func{\ell(\CC^*;p)\times \ell(\CC;p)}{\CC}$ given by $Q_p(a;b)=q_{p,1}(a;b;a)$, then $Q_p$  is determined by the following non-linear recurrence
\begin{equation}
\label{p=1,2}
Q_1(a;b)=\dfrac{b(0)}{2a(0)}, \hspace{.5cm}Q_2(a;b)=\dfrac{1}{2a(0)a(1)}\Big(b(0)b(1)- a(0)^2-a(1)^2\Big), 
\end{equation}
and
\begin{equation}
\label{recurrence}
Q_{p+1}(a;b)=
\left\{\begin{array}{cl}
Q_{p}(\hat a;\hat b), &  b(p)\not=0,\\[1ex]
Q_{p}(\check a;\check b), &   b(0)\not=0,\\[1ex]
Q_{p-1}(\tilde a;\tilde b), &  b(p)=b(0)=0,
\end{array}\right.
\end{equation}
for $p\ge 2$, where the periodic sequences $\hat a,\hat b,\check a,\check b,\tilde a$ and $\tilde b$ are defined as  the periodic extension of 
$$
\begin{array}{llll}
\hat a(k)=a(k),& \hspace{-.15cm}k=0,\ldots,p-2, &  \hat a(p-1)=\dfrac{a(p-1)a( p)}{b( p)},& \\[3ex]
\hat b(k)=b(k), &\hspace{-.15cm} k=1,\ldots,p-2,  &  \hat b(p-1)=\dfrac{b(p-1)b( p)-a(p-1)^2}{b( p)},&  \hat b(0)=\dfrac{b(0)b( p)-a( p)^2}{b( p)},\\[3ex]
\check a(k)=a(k+1),&\hspace{-.15cm} k=0,\ldots,p-2, &  \check a(p-1)=\dfrac{a(0)a( p)}{b(0)},& \\[3ex]
\check b(k)=b(k+1), & \hspace{-.15cm}k=1,\ldots,p-2,  &  \check b(p-1)=\dfrac{b(0)b( p)-a( p)^2}{b(0)},&  \check b(0)=\dfrac{b(0)b(1)-a(0)^2}{b(0)},\\[3ex]
\tilde a(k)=a(k+1),&\hspace{-.15cm} k=0,\ldots,p-3, &  \tilde a(p-2)=-\dfrac{a(0)a(p-1)}{a( p)},& 
\\[3ex]
\tilde b(k)=b(k+1), & \hspace{-.15cm}k=0,\ldots,p-2.  & &
\end{array}
$$
Notice that when $p=2$, $\tilde a(0)=-\dfrac{a(0)a(1)}{a( 2)}$.
In addition, the properties of Chebyshev equations and its solutions, establishes that when $b(0),b(p)\not=0$, necessarily $Q_{p}(\hat a;\hat b)=Q_{p}(\check a;\check b)$. \vspace{1cc}

Our next aim is to   obtain a closed expression for the  Floquet functions and with this end, we  introduce some concepts and notations. 

A {\it binary multi-index of order $p$} is a $p$-tuple $\alpha=(\alpha_0,\ldots,\alpha_{p-1})\in \{0,1\}^p$ and its {\it length} is defined as  $|\alpha|=\sum\limits_{j=0}^{p-1}\alpha_j\le p$. So $|\alpha|=m$ iff   exactly $m$ components of $\alpha$ are equal to $1$ and exactly $p-m$ components of $\alpha$ are equal to $0$. The only binary multi-index of order $p$ whose length equals $p$ is $\pi_p=(1,\ldots,1)$. Moreover, $\alpha\in \{0,1\}^p$ and $|\alpha|=m$, we consider $0\le i_1<\cdots<i_m\le p-1$ such that $\alpha_{i_1}=\cdots=\alpha_{i_m}=1$. 

Given $\alpha\in \{0,1\}^p$ and a sequence $a\in \ell(\KK)$, we consider the following values
\begin{equation}
\label{exp}
a^\alpha=\prod\limits_{j=0}^{p-1}a(j)^{\alpha_j}\hspace{.5cm}\hbox{and}\hspace{.5cm}a^{2\alpha}=\prod\limits_{j=0}^{p-1}a(j)^{2\alpha_j}
\end{equation}
respectively, where we assume  $0^0=1$. Observe that $a^{\pi_p}=\prod\limits_{j=0}^{p-1}a(j)$ for any $a\in \ell(\KK;p)$. In the sequel we also assume   the usual convention that empty sums and empty products are defined as $0$ and $1$, respectively.

Given $p\in \NN^*$, we define $\Lambda_p^0=\{(0,\ldots,0)\}$ and for $p\ge 2$, 
$\Lambda_p^1=\big\{\alpha\in \{0,1\}^p: |\alpha|=1\big\}$. Moreover, when $p\ge 4$, for any $m=2,\ldots,\lfloor\frac{p}{2}\rfloor$, we define
\begin{equation}
\label{floor}
\Lambda_p^m=\Big\{\alpha\in \{0,1\}^p: |\alpha|=m,\hspace{.15cm}\hbox{and $i_j+2\le i_{j+1}\le p-2(m-j)+\min\{i_1,1\}$, \hspace{.15cm}$j=1,\ldots,m-1$}\Big\},
\end{equation}
or, equivalently,  
$$\Lambda_p^m=\Big\{\alpha\in \{0,1\}^p: |\alpha|=m,\hspace{.15cm}\hbox{ $i_{j+1}-i_j\ge 2$, \hspace{.15cm}$j=1,\ldots,m-1$ and $i_m\le p-2$ when $i_1=0$}\Big\}.
$$
\vspace{.5cc}

Given $p\ge 2$, $m=1,\ldots,\lfloor\frac{p}{2}\rfloor$ and $\alpha\in \Lambda_p^m$, let $0\le i_1<\cdots<i_m\le p-1$ be the indexes such that  $\alpha_{i_1}=\cdots=\alpha_{i_m}=1$. Then, we define the binary multi-index $\bar \alpha$ of order $p$ as 
$$\bar \alpha_{i_j}=\bar \alpha_{i_j+1}=0, \hspace{.15cm}j=1,\ldots,m,\hspace{.25cm}\hbox{and}\hspace{.25cm}\bar \alpha_i=1\hspace{.25cm}\hbox{otherwise},$$
where if $i_m=p-1$, then $\bar \alpha_{p-1}=\bar \alpha_0=0$.   Moreover, if $\alpha\in\Lambda_p^0$; that is, if $\alpha=(0,\ldots,0)$, then we define $\bar \alpha=\pi_p$. It is clear that, in any case,   $|\bar \alpha|=p-2m$.

Given $p\ge 2$ and $0\le j\le\lfloor\frac{p}{2}\rfloor$, we define the following sets of binary multi-indices of order $p$
$$\begin{array}{rlrl}
A^{j,1}_p=&\hspace{-.25cm}\big\{\alpha \in \Lambda_{p}^{j}: \bar \alpha_0=\bar\alpha_{p-1}=1\big\},&\\[1ex]
 A^{j,2}_p=&\hspace{-.25cm} \big\{\alpha \in \Lambda_{p}^{j}: \bar \alpha_0=0,\,\,\bar\alpha_{p-1}=1\big\},&  A^{j,3}_p=&\hspace{-.25cm}\big\{\alpha \in \Lambda_{p}^{j}: \bar \alpha_0=1,\,\,\bar\alpha_{p-1}=0\big\}, \\[1ex] 
 A^{j,4}_p=&\hspace{-.25cm}\big\{\alpha \in \Lambda_{p}^{j}: \bar \alpha_0=\bar\alpha_{p-1}=0\hspace{.15cm}\hbox{and}\hspace{.15cm}\alpha_0=1\big\}, & A^{j,5}_p=&\hspace{-.25cm}\big\{\alpha \in \Lambda_{p}^{j}: \bar \alpha_0=\bar\alpha_{p-1}=0\hspace{.15cm}\hbox{and}\hspace{.15cm}\alpha_0=0\big\},\end{array}$$
 that clearly determine a partition of $\Lambda_p^j$.  Moreover, $A^{0,1}_p=\{(0,\ldots,0)\}$ and  $A^{0,2}_{p}=A^{0,3}_{p}=A^{0,4}_{p}=A^{0,5}_{p}=\emptyset$.
%
 \begin{lemma}
\label{partition}
Given $p\ge 4$ and $2\le j\le\lfloor\frac{p}{2}\rfloor$, then $A^{j,3}_{p+1}=A^{j,5}_p\times \{0\}$ and moreover 
$$\begin{array}{rlrl}
A^{j,1}_{p+1}=&\hspace{-.25cm}\big(A^{j,1}_p\times \{0\}\big)\cup \big(A^{j,3}_p\times \{0\}\big), & \hspace{.5cm}
 A^{j,2}_{p+1}=&\hspace{-.25cm} \big(A^{j,2}_p\times \{0\}\big)\cup \big(A^{j,4}_p\times \{0\}\big), \\[1ex]
 A^{j,4}_{p+1}=&\hspace{-.25cm}\big(A^{j-1,2}_{p-1}\times \{(1,0)\}\big)\cup \big(A^{j-1,4}_{p-1}\times \{(1,0)\}\big), & \hspace{.5cm}
 A^{j,5}_{p+1}=&\hspace{-.25cm}\big(A^{j-1,1}_p\times \{1\}\big)\cup \big(A^{j-1,3}_p\times \{1\}\big).\end{array}$$
\end{lemma}

\begin{propo}
\label{cardinal}
Given $p\in \NN^*$ and $0\le j\le\lfloor\frac{p}{2}\rfloor$, then $\displaystyle |\Lambda^{j}_p|=\frac{p}{p-j}{p-j\choose j}$.  Therefore, for any $m\in \NN^*$ we have that $\displaystyle\sum\limits_{j=0}^m|\Lambda^{j}_{2m}|=2T_{m}\Big(\frac{3}{2}\Big)$ and $\displaystyle \sum\limits_{j=0}^m|\Lambda^{j}_{2m+1}|=W_{m}\Big(\frac{3}{2}\Big)$. 
\end{propo}
\proof We know that  $|\Lambda_p^0|=1$, for any $p\in \NN^*$ and that $|\Lambda_p^1|=\big|\big\{\alpha\in \{0,1\}^p: |\alpha|=1\big\}\big|=p$, for any $p\ge 2$.

 If $\alpha\in \Lambda_{2m}^m$, $m\in \NN^*$, and $0\le i_1<\cdots<i_m\le 2m-1$ are such that  $\alpha_{i_1}=\cdots=\alpha_{i_m}=1$, then  $0\le i_1\le 1$ and 
$$2+i_j\le i_{j+1}\le 2j+\min\{i_1,1\},\hspace{.25cm}j=1,\ldots,m-1.$$
 
If $i_1=0$, then $i_j=2(j-1)$, $j=1,\ldots,m$, whereas  when $i_1=1$, then 
$i_j=2(j-1)+1$, $j=1,\ldots,m$. In both cases $\bar \alpha=(0,\ldots,0)$ and moreover $|\Lambda_{2m}^m|=2$.

If $\alpha\in \Lambda_{2m+1}^m$, then  $0\le i_1<\cdots<i_m\le 2m$ are such that  $\alpha_{i_1}=\cdots=\alpha_{i_m}=1$, then  $0\le i_1\le 2$ and 
$$2+i_j\le i_{j+1}\le 2j+1+\min\{i_1,1\},\hspace{.25cm}j=1,\ldots,m-1.$$
 
If $i_1=0$, then either $i_j=2(j-1)$, $j=1,\ldots,m$, which implies that $\bar \alpha=(0,\ldots,0,1)$; or there exists $2\le \ell\le m$ such that $i_j=2(j-1)$ when $1\le j< \ell$ and $i_\ell=2\ell-1$. Then $i_j=2j-1$, $j=\ell,\ldots,m$ and hence, $\bar \alpha_{2\ell-2}=1$ and $\bar \alpha_i=0$, otherwise. Then, $|\{\alpha \in \Lambda_{2m+1}^m:\alpha_0=1\}|=m$.

If $i_1=1$, then either $i_j=2j-1$, $j=1,\ldots,m$, which implies that $\bar \alpha=(1,0\ldots,0)$; or there exists $2\le \ell\le m$ such that $i_j=2j-1$ when $1\le j< \ell$ and $i_\ell=2\ell$. Then $i_j=2j$, $j=\ell,\ldots,m$ and hence, $\bar \alpha_{2\ell-1}=1$ and $\bar \alpha_i=0$, otherwise. Moreover, $|\{\alpha \in \Lambda_{2m+1}^m:\alpha_0=0,\,\alpha_1=1\}|=1$.

If $i_1=2$, then $i_j=2j$, $j=1,\ldots,m$ and hence, $\bar \alpha_{1}=1$ and $\bar \alpha_i=0$, otherwise; which in turns implies that $|\{\alpha \in \Lambda_{2m+1}^m:\alpha_0=\alpha_1=0,\,\alpha_2=1\}|=1$. 

Therefore,  $|\{\alpha \in \Lambda_{2m+1}^m:\alpha_1=1\}|=2m+1$ and hence, we have obtained that given $p\in \NN^*$, the claimed formula for $|\Lambda^j_p|$ is true for $j=0,1$ and for $j=\lfloor\frac{p}{2}\rfloor$. 

Assume now that the formula is true for $p\ge 2$ and  $0\le j\le\lfloor\frac{p}{2}\rfloor$. Then given $1\le j\le\lfloor\frac{p+1}{2}\rfloor-1$ and applying Lemma \ref{partition}, we get that 
$$\begin{array}{rl}
|\Lambda^{j+1}_{p+1}|=&\hspace{-.25cm}\displaystyle \sum\limits_{i=1}^3|A^{j+1,i}_{p+1}|+\sum\limits_{i=4}^5|A^{j+1,i}_{p+1}|=\sum\limits_{i=1}^5|A^{j+1,i}_{p}|+   |A^{j,2}_{p-1}|+|A^{j,4}_{p-1}|+ |A^{j,1}_{p}|+ |A^{j,3}_{p}|\\[3ex]
=&\hspace{-.25cm}\displaystyle|\Lambda^{j+1}_{p}|+ |A^{j,2}_{p-1}|+|A^{j,4}_{p-1}|+ |A^{j,1}_{p-1}|+ |A^{j,3}_{p-1}|+ |A^{j,5}_{p-1}|=|\Lambda^{j+1}_{p}|+|\Lambda^{j}_{p-1}|\\[2ex]
=&\hspace{-.25cm}\displaystyle\frac{p}{p-j-1}{p -j-1\choose j+1}+\frac{p-1}{p-1-j}{p-1-j\choose j}=\dfrac{p+1}{p-j}{p-j\choose j+1}.\end{array}$$

Finally, from the Identity \eqref{ch:first} we obtain that 
$T_p\big(\frac{i}{2}\big)=\dfrac{(i)^p}{2}\sum\limits_{j=0}^{\lfloor \frac{p}{2}\rfloor}|\Lambda^{j}_{p}|$   for any $p\in \NN^*$,  that  in turns,  from the identities \eqref{doubling:first} and taking into account that $T_m(-x)=(-1)^mT_m(x)$ and $V_m(-x)=(-1)^mW_m(x)$  for any $m\in \NN^*$, implies the last claims. \qed 
\vspace{.5cc}

\begin{corollary}
\label{even:odd} Given $p,m\in \NN^*$, then the following identities hold
$$\begin{array}{rl}
\displaystyle \sum\limits_{\alpha\in \Lambda^0_{p}}a^{2\alpha}b^{\bar \alpha}=&\hspace{-.25cm}\displaystyle b^{\pi_p}=\prod\limits_{j=0}^{p-1}b(j)\\[3ex]
\displaystyle \sum\limits_{\alpha\in \Lambda^1_{p}}a^{2\alpha}b^{\bar \alpha}=&\hspace{-.25cm}\displaystyle 
\sum\limits_{i=0}^{p-2}a(i)^2\prod\limits_{j=0\atop j\not=i,i+1}^{p-1}b(j)+a(p-1)^2\prod\limits_{j=1}^{p-2}b(j)\\[3ex]
\displaystyle \sum\limits_{\alpha\in \Lambda^m_{2m}}a^{2\alpha}b^{\bar \alpha}=&\hspace{-.25cm}\displaystyle \prod\limits_{j=0}^{m-1}a(2j)^2+\prod\limits_{j=1}^{m}a(2j-1)^2\\[3ex]
\displaystyle \sum\limits_{\alpha\in \Lambda^m_{2m+1}}a^{2\alpha}b^{\bar \alpha}=&\hspace{-.25cm}\displaystyle \sum\limits_{i=0}^{m}b(2i)\prod\limits_{j=0}^{i-1}a(2j)^2\prod\limits_{j=i+1}^{m}a(2j-1)^2+\sum\limits_{i=1}^{m}b(2i-1)\prod\limits_{j=1}^{i-1}a(2j-1)^2\prod\limits_{j=i}^{m}a(2j)^2\end{array}$$
\end{corollary}
\vspace{.5cc}

\begin{theorem}
\label{Ff} 
Given $p\in \NN^*$ and  $r\in \CC^*$ then for any $a,c\in \ell(\CC^*;p,r)$ and  $b\in \ell(\CC;p,r)$,
we have that 
$$q_{p,r}(a;b;c)=\dfrac{1}{2 }\sqrt{\dfrac{r}{a^{\pi_p} c^{\pi_p}}}
\sum\limits_{j=0}^{\lfloor\frac{p}{2}\rfloor}(-1)^j\sum\limits_{\alpha\in \Lambda_p^j}r^{-\alpha_{p-1}}a^{\alpha} b^{\bar \alpha} c^\alpha.$$\end{theorem}
\proof We first prove, by induction on $p$, that for any $a\in \ell(\CC^*,p)$ and $b\in \ell(\CC;p)$ we have 
$$Q_p(a;b)=q_{p,1}(a;b;a)=\dfrac{1}{2a^{\pi_p}}
\sum\limits_{j=0}^{\lfloor\frac{p}{2}\rfloor}(-1)^j\sum\limits_{\alpha\in \Lambda_p^j}a^{2\alpha} b^{\bar \alpha}.$$
 
From Corollary  \ref{even:odd}, for $p=1$ the claimed formula gives the value $\dfrac{b(0)}{2a(0)}$, whereas for $p=2$ gives the value $\dfrac{1}{2a(0)a(1)}\Big(b(0)b(1)-a(0)^2-a(1)^2\Big)$. Therefore, taking into account  the identities \eqref{p=1,2}, the proposed formula coincides with the expression for $Q_p$ where $p=1,2$. Assume now that it is true for $p\ge 2$ and consider $a\in \ell(\KK^*;p+1)$ and $b\in \ell(\KK;p+1)$. Since the hypotheses $b(0)\not=0$ or  $b(0)=b(p)=0$ lead  to analogous reasoning that the case $b(p)\not=0$, in the sequel we always assume that  $b(p)\not=0$ and hence, our aim is to prove that $Q_p(\hat a;\hat b)=Q_{p+1}(a;b)$ for any $p\ge 2$. 

When $p=2$, then 
$$\begin{array}{rl}
Q_2(\hat a;\hat b)=&\hspace{-.25cm}\dfrac{b(2)}{2a(0)a(1)a(2)}\left(\dfrac{\big(b(0)b(2)-a(2)^2)\big(b(1)b(2)-a(1)^2\big)}{b(2)^2}-a(0)^2-\dfrac{a(1)^2a(2)^2}{b(2)^2}\right)\\[3ex]
=&\hspace{-.25cm}\dfrac{1}{2a(0)a(1)a(2)}\Big(b(0)b(1)b(2)-a(1)^2b(0)-a(2)^2b(1)-a(0)^2b(2)\Big)=Q_3(a;b).\end{array}$$

When $p=3$, then 
$$\begin{array}{rl}
Q_3(\hat a;\hat b)=&\hspace{-.25cm}\dfrac{b(3)}{2a(0)a(1)a(2)a(3)}\left(\dfrac{b(1)\big(b(0)b(3)-a(3)^2)\big(b(2)b(3)-a(2)^2\big)}{b(3)^2}-\dfrac{a(1)^2\big(b(0)b(3)-a(3)^2)}{b(3)}\right.\\[3ex]
&\hspace{3cm}\left.-\dfrac{a(2)^2a(3)^2b(1)}{b(3)^2}-\dfrac{a(0)^2\big(b(2)b(3)-a(2)^2)}{b(3)}\right)\\[3ex]
=&\hspace{-.25cm}\dfrac{1}{2a(0)a(1)a(2)a(3)}\Big(b(0)b(1)b(2)b(3)-a(2)^2b(0)b(1)-a(3)^2b(1)b(2)-a(1)^2b(0)b(3)\\[3ex]
& \hspace{3cm} +a(1)^2a(3)^2-a(0)^2b(2)b(3)+a(0)^2a(2)^2\Big)=Q_4(a;b).\end{array}$$

When $p\ge 4$,   taking into account that 
$$b(p) \hat a^{\pi_p}= a^{\pi_{p+1}}\hspace{.25cm}\hbox{and}\hspace{.25cm}b(p)\hat b^{\pi_p}=b^{\pi_{ p+1}}- a(p-1)^2\prod\limits_{j=0}^{p-2}b(j)-a(p)^2\prod\limits_{j=1}^{p-1}b(j)+\dfrac{a(p-1)^2a(p)^2}{b(p)}\prod\limits_{j=1}^{p-2}b(j),$$ 
the first two identities of Lemma \ref{even:odd} imply that 
$$ Q_p(\hat a;\hat b)=\dfrac{1}{2a^{\pi_{p+1}}}\left[ \sum\limits_{\alpha \in \Lambda_{p+1}^0}a^{2\alpha}b^{\bar \alpha}-\sum\limits_{\alpha \in \Lambda_{p+1}^1}a^{2\alpha}b^{\bar \alpha}+b(p)\sum\limits_{j=2}^{\lfloor\frac{p}{2}\rfloor}(-1)^j\sum\limits_{\alpha\in \Lambda_p^j}\hat a^{2\alpha} \hat b^{\bar \alpha} +R_p(a;b) \right]$$
%
where 
$$\begin{array}{rl}
R_p(a;b)=&\hspace{-.25cm}\displaystyle a(p-1)^2\sum\limits_{i=0}^{p-3} a(i)^2\prod\limits_{j=0\atop j\not=i,i+1}^{p-2} b(j)+a(p)^2\sum\limits_{i=1}^{p-2} a(i)^2\prod\limits_{j=1\atop j\not=i,i+1}^{p-1} b(j)-\dfrac{a(p-1)^2a(p)^2}{b(p)}\sum\limits_{i=1}^{p-3} a(i)^2\!\!\!\prod\limits_{j=1\atop j\not=i,i+1}^{p-2} b(j)\\[3ex]
=&\hspace{-.25cm}\displaystyle \sum\limits_{\alpha \in A^{2,3}_{p+1}}a^{2\alpha}b^{\bar \alpha}+\sum\limits_{\alpha \in A^{2,4}_{p+1}}a^{2\alpha}b^{\bar \alpha}+\sum\limits_{\alpha \in A^{2,5}_{p+1}}a^{2\alpha}b^{\bar \alpha}-\dfrac{a(p)^2}{b(p)}\sum\limits_{\alpha \in A^{2,5}_p}a^{2\alpha}b^{\bar \alpha}.\end{array}$$ 

On the other hand, from the two last identities of Corollary \ref{even:odd}, when $p$ is even then 
$$\begin{array}{rl}
\displaystyle b(p)\sum\limits_{\alpha\in \Lambda^{\lfloor\frac{p}{2}\rfloor}_{p}}\hat a^{2\alpha}\hat b^{\bar \alpha}=&\hspace{-.25cm}\displaystyle b(p)\prod\limits_{j=0}^{\lfloor\frac{p}{2}\rfloor-1}a(2j)^2+\dfrac{a(p)^2}{b(p)}\prod\limits_{j=1}^{\lfloor\frac{p}{2}\rfloor}a(2j-1)^2\\[3ex]
=&\hspace{-.25cm}\displaystyle \sum\limits_{\alpha \in A^{\lfloor\frac{p}{2}\rfloor,1}_{p+1}}a^{2\alpha}b^{\bar \alpha}+\sum\limits_{\alpha \in A^{\lfloor\frac{p}{2}\rfloor,2}_{p+1}}a^{2\alpha}b^{\bar \alpha}+\dfrac{a(p)^2}{b(p)}\sum\limits_{\alpha \in A^{\lfloor\frac{p}{2}\rfloor,5}_{p}}a^{2\alpha}b^{\bar \alpha}\end{array}$$
since $A^{\lfloor\frac{p}{2}\rfloor,1}_{p+1}=\emptyset$; whereas when $p$ is odd, then $\lfloor\frac{p+1}{2}\rfloor=\lfloor\frac{p}{2}\rfloor+1$ and 
$$\begin{array}{rl}
\displaystyle b(p)\sum\limits_{\alpha\in \Lambda^{\lfloor\frac{p}{2}\rfloor}_{p}}\hat a^{2\alpha}\hat b^{\bar \alpha}=&\hspace{-.25cm}\displaystyle 
\big(b(0)b(p)-a(p)^2\big)\prod\limits_{j=1}^{\lfloor\frac{p}{2}\rfloor}a(2j-1)^2+\big(b(p-1)b(p)-a(p-1)^2\big)\prod\limits_{j=0}^{\lfloor\frac{p}{2}\rfloor-1} a(2j)^2\\[3ex]
+&\hspace{-.25cm}\displaystyle  b(p)\sum\limits_{i=1}^{\lfloor\frac{p}{2}\rfloor-1} b(2i)\prod\limits_{j=0}^{i-1} a(2j)^2\prod\limits_{j=i+1}^{\lfloor\frac{p}{2}\rfloor}a(2j-1)^2+\dfrac{a(p)^2}{b(p)}\sum\limits_{i=1}^{\lfloor\frac{p}{2}\rfloor} b(2i-1)\prod\limits_{j=1}^{i-1} a(2j-1)^2\prod\limits_{j=i}^{\lfloor\frac{p}{2}\rfloor} a(2j)^2\\[3ex] 
=&\hspace{-.25cm}\displaystyle 
 -\prod\limits_{j=1}^{\lfloor\frac{p+1}{2}\rfloor}a(2j-1)^2-\prod\limits_{j=0}^{\lfloor\frac{p}{2}\rfloor} a(2j)^2\\[3ex]
+&\hspace{-.25cm}\displaystyle  b(p)\sum\limits_{i=0}^{\lfloor\frac{p}{2}\rfloor} b(2i)\prod\limits_{j=0}^{i-1} a(2j)^2\prod\limits_{j=i+1}^{\lfloor\frac{p}{2}\rfloor}a(2j-1)^2+\dfrac{a(p)^2}{b(p)}\sum\limits_{i=1}^{\lfloor\frac{p}{2}\rfloor} b(2i-1)\prod\limits_{j=1}^{i-1} a(2j-1)^2\prod\limits_{j=i}^{\lfloor\frac{p}{2}\rfloor} a(2j)^2\\[3ex] 
=&\hspace{-.25cm}\displaystyle 
 -\sum\limits_{\alpha\in \Lambda^{\lfloor\frac{p+1}{2}\rfloor}_{p+1}}a^{2\alpha}b^{\bar \alpha}+\sum\limits_{\alpha\in A^{\lfloor\frac{p}{2}\rfloor,1}_{p+1}}a^{2\alpha}b^{\bar \alpha}+\sum\limits_{\alpha\in A^{\lfloor\frac{p}{2}\rfloor,2}_{p+1}}a^{2\alpha}b^{\bar \alpha}+\dfrac{a(p)^2}{b(p)}\sum\limits_{\alpha \in A^{\lfloor\frac{p}{2}\rfloor,5}_p}a^{2\alpha}b^{\bar \alpha}.\end{array}$$

In particular, when $p=4,5$ then $\lfloor\frac{p}{2}\rfloor=2$  we obtain that 
$$b(4)\sum\limits_{\alpha\in \Lambda_4^2}\hat a^{2\alpha} \hat b^{\bar \alpha} +R_4(a;b) =\sum\limits_{\alpha\in \Lambda_{5}^2}a^{2\alpha}  b^{\bar \alpha}\hspace{.25cm}\hbox{and}\hspace{.25cm}b(5)\sum\limits_{\alpha\in \Lambda_5^2}\hat a^{2\alpha} \hat b^{\bar \alpha} +R_5(a;b) =-\sum\limits_{\alpha\in \Lambda_{6}^3}a^{2\alpha}  b^{\bar \alpha}+\sum\limits_{\alpha\in \Lambda_{6}^2}a^{2\alpha}  b^{\bar \alpha}$$
and hence $Q_p(\hat a;\hat b)=Q_{p+1}( a;b)$.

Consider now $p\ge 6$ and $2\le j\le \lfloor\frac{p}{2}\rfloor-1$. Then, 
$b(p)\sum\limits_{\alpha \in \Lambda^j_p}\hat a^{2\alpha}\hat b^{\bar \alpha}=b(p)\sum\limits_{i=1}^5\sum\limits_{\alpha \in A^{j,i}_{p}}\hat a^{2\alpha}\hat b^{\bar \alpha}$ and moreover we have
$$\begin{array}{rl}
\displaystyle b(p)\sum\limits_{\alpha \in A^{j,1}_p}\hat a^{2\alpha}\hat b^{\bar \alpha}=&\hspace{-.25cm}\displaystyle  
b(p)\sum\limits_{\alpha \in A^{j,1}_p}  a^{2\alpha}b^
{\bar \alpha}-a(p-1)^2 \sum\limits_{\alpha \in A^{j,1}_p}  a^{2\alpha}\prod\limits_{i=0}^{p-2} b^{\bar \alpha_i}-a(p)^2 \sum\limits_{\alpha \in A^{j,1}_p}  a^{2\alpha}\prod\limits_{i=1}^{p-1} b^{\bar \alpha_i}\\[4ex]
+ &\hspace{-.25cm}\displaystyle  \dfrac{a(p-1)^2a(p)^2}{b(p)}\sum\limits_{\alpha \in A^{j,1}_p}  a^{2\alpha}\prod\limits_{i=1}^{p-2} b^{\bar \alpha_i},\\[4ex]
=&\hspace{-.25cm}\displaystyle  
\sum\limits_{\alpha \in A^{j,1}_p\times\{0\}}  a^{2\alpha}b^
{\bar \alpha}-\sum\limits_{\alpha \in A^{j+1,5}_{p}\times \{0\}}  a^{2\alpha} b^{\bar \alpha}- \sum\limits_{\alpha \in A^{j,1}_p\times \{1\}}  a^{2\alpha}b^{\bar \alpha}+\dfrac{a(p)^2}{b(p)}\sum\limits_{\alpha \in A^{j+1,5}_p}  a^{2\alpha}b^{\bar \alpha},\\[4ex]
\displaystyle b(p)\sum\limits_{\alpha \in A^{j,2}_p}\hat a^{2\alpha}\hat b^{\bar \alpha}=&\hspace{-.25cm}\displaystyle  
b(p)\sum\limits_{\alpha \in A^{j,2}_p}  a^{2\alpha}b^
{\bar \alpha}-a(p-1)^2 \sum\limits_{\alpha \in A^{j,2}_p}  a^{2\alpha}\prod\limits_{i=0}^{p-2} b^{\bar \alpha_i}=\sum\limits_{\alpha \in A^{j,2}_p\times \{0\}}  a^{2\alpha}b^
{\bar \alpha}- \sum\limits_{\alpha \in A^{j+1,4}_{p+1}}  a^{2\alpha}  b^{\bar \alpha },\\[4ex]
\displaystyle b(p)\sum\limits_{\alpha \in A^{j,3}_p}\hat a^{2\alpha}\hat b^{\bar \alpha}=&\hspace{-.25cm}\displaystyle  
b(p)\sum\limits_{\alpha \in A^{j,3}_p}  a^{2\alpha}b^
{\bar \alpha}-a(p)^2 \sum\limits_{\alpha \in A^{j,3}_p}  a^{2\alpha}\prod\limits_{i=1}^{p-1} b^{\bar \alpha_i}=\sum\limits_{\alpha \in A^{j,3}_p\times \{0\}}  a^{2\alpha}b^
{\bar \alpha}-\sum\limits_{\alpha \in A^{j,3}_p\times \{1\}}  a^{2\alpha}b^{\bar \alpha},\\[4ex]
\displaystyle b(p)\sum\limits_{\alpha \in A^{j,4}_p}\hat a^{2\alpha}\hat b^{\bar \alpha}=&\hspace{-.25cm}\displaystyle b(p)\sum\limits_{\alpha \in A^{j,4}_p} a^{2\alpha} b^{\bar \alpha}=\sum\limits_{\alpha \in A^{j,4}_p\times \{0\}} a^{2\alpha} b^{\bar \alpha},\\[4ex]
\displaystyle b(p)\sum\limits_{\alpha \in A^{j,5}_p}\hat a^{2\alpha}\hat b^{\bar \alpha}=&\hspace{-.25cm}\displaystyle  \dfrac{a(p)^2}{b(p)}\sum\limits_{\alpha \in A^{j,5}_p} a^{2\alpha} b^{\bar \alpha}\end{array}$$
Therefore, from Lemma \ref{partition}  we obtain that 
$$\begin{array}{rl}
\displaystyle b(p)\sum\limits_{j=2}^{\lfloor\frac{p}{2}\rfloor-1}(-1)^j\sum\limits_{\alpha\in \Lambda_p^j}\hat a^{2\alpha} \hat b^{\bar \alpha} =&\hspace{-.25cm}\displaystyle \sum\limits_{j=2}^{\lfloor\frac{p}{2}\rfloor-1}(-1)^j\Bigg[\sum\limits_{\alpha \in A^{j,1}_{p+1}}  a^{2\alpha}b^
{\bar \alpha}+\sum\limits_{\alpha \in A^{j,2}_{p+1}}  a^{2\alpha}b^
{\bar \alpha}\Bigg]\\[4ex]
-&\hspace{-.25cm}\displaystyle\sum\limits_{j=2}^{\lfloor\frac{p}{2}\rfloor-1}(-1)^j\Bigg[\sum\limits_{\alpha \in A^{j+1,3}_{p+1}}  a^{2\alpha}b^
{\bar \alpha}+\sum\limits_{\alpha \in A^{j+1,4}_{p+1}}  a^{2\alpha}b^
{\bar \alpha}+\sum\limits_{\alpha \in A^{j+1,5}_{p+1}}  a^{2\alpha}b^
{\bar \alpha}\Bigg]\\[4ex]
+&\hspace{-.25cm}\displaystyle \dfrac{a(p)^2}{b(p)}\sum\limits_{j=2}^{\lfloor\frac{p}{2}\rfloor-1}(-1)^j\Bigg[\sum\limits_{\alpha \in A^{j+1,5}_{p}}a^{2\alpha}b^{\bar \alpha}+
\sum\limits_{\alpha \in A^{j,5}_{p}}a^{2\alpha}b^{\bar \alpha}\Bigg]\\[4ex]
=&\hspace{-.25cm}\displaystyle \sum\limits_{\alpha \in A^{2,1}_{p+1}}  a^{2\alpha}b^
{\bar \alpha}+\sum\limits_{\alpha \in A^{2,2}_{p+1}}  a^{2\alpha}b^
{\bar \alpha}+\sum\limits_{j=3}^{\lfloor\frac{p}{2}\rfloor-1}(-1)^j \sum\limits_{\alpha \in \Lambda^{j}_{p+1}}  a^{2\alpha}b^
{\bar \alpha}\\[4ex]
+&\hspace{-.25cm}\displaystyle (-1)^{\lfloor\frac{p}{2}\rfloor}\Bigg[\sum\limits_{\alpha \in A^{{\lfloor\frac{p}{2}\rfloor},3}_{p+1}}  a^{2\alpha}b^
{\bar \alpha}+\sum\limits_{\alpha \in A^{{\lfloor\frac{p}{2}\rfloor},4}_{p+1}}  a^{2\alpha}b^
{\bar \alpha}+\sum\limits_{\alpha \in A^{{\lfloor\frac{p}{2}\rfloor},5}_{p+1}}  a^{2\alpha}b^
{\bar \alpha}\Bigg]\\[4ex]
+&\hspace{-.25cm}\displaystyle \dfrac{a(p)^2}{b(p)}\Bigg[
\sum\limits_{\alpha \in A^{2,5}_{p}}a^{2\alpha}b^{\bar \alpha}-(-1)^{\lfloor\frac{p}{2}\rfloor}\sum\limits_{\alpha \in A^{\lfloor\frac{p}{2}\rfloor,5}_{p}}a^{2\alpha}b^{\bar \alpha}\Bigg],
\end{array}$$
$$\begin{array}{rl}
R_p(a;b)=&\hspace{-.25cm}\displaystyle a(p-1)^2\sum\limits_{i=0}^{p-3} a(i)^2\prod\limits_{j=0\atop j\not=i,i+1}^{p-2} b(j)+a(p)^2\sum\limits_{i=1}^{p-2} a(i)^2\prod\limits_{j=1\atop j\not=i,i+1}^{p-1} b(j)-\dfrac{a(p-1)^2a(p)^2}{b(p)}\sum\limits_{i=1}^{p-3} a(i)^2\!\!\!\prod\limits_{j=1\atop j\not=i,i+1}^{p-2} b(j)\\[3ex]
=&\hspace{-.25cm}\displaystyle \sum\limits_{\alpha \in A^{2,3}_{p+1}}a^{2\alpha}b^{\bar \alpha}+\sum\limits_{\alpha \in A^{2,4}_{p+1}}a^{2\alpha}b^{\bar \alpha}+\sum\limits_{\alpha \in A^{2,5}_{p+1}}a^{2\alpha}b^{\bar \alpha}-\dfrac{a(p)^2}{b(p)}\sum\limits_{\alpha \in A^{2,5}_p}a^{2\alpha}b^{\bar \alpha},\end{array}$$ 
which implies that 
$$\begin{array}{rl}
\displaystyle b(p)\sum\limits_{j=2}^{\lfloor\frac{p}{2}\rfloor}(-1)^j\sum\limits_{\alpha\in \Lambda_p^j}\hat a^{2\alpha} \hat b^{\bar \alpha} +R_p(a;b)=&\hspace{-.25cm}\displaystyle \sum\limits_{j=2}^{\lfloor\frac{p}{2}\rfloor-1}(-1)^j \sum\limits_{\alpha \in \Lambda^{j}_{p+1}}  a^{2\alpha}b^
{\bar \alpha}\\[4ex]
+&\hspace{-.25cm}\displaystyle (-1)^{\lfloor\frac{p}{2}\rfloor}\Bigg[\sum\limits_{\alpha \in A^{{\lfloor\frac{p}{2}\rfloor},3}_{p+1}}  a^{2\alpha}b^
{\bar \alpha}+\sum\limits_{\alpha \in A^{{\lfloor\frac{p}{2}\rfloor},4}_{p+1}}  a^{2\alpha}b^
{\bar \alpha}+\sum\limits_{\alpha \in A^{{\lfloor\frac{p}{2}\rfloor},5}_{p+1}}  a^{2\alpha}b^
{\bar \alpha}\Bigg]\\[4ex]
+&\hspace{-.25cm}\displaystyle (-1)^{\lfloor\frac{p}{2}\rfloor}\Bigg[b(p)\sum\limits_{\alpha\in \Lambda^{\lfloor\frac{p}{2}\rfloor}_{p}}\hat a^{2\alpha}\hat b^{\bar \alpha}-\dfrac{a(p)^2}{b(p)}\sum\limits_{\alpha \in A^{\lfloor\frac{p}{2}\rfloor,5}_{p}}a^{2\alpha}b^{\bar \alpha}\Bigg].
\end{array}$$

Finally, using the expression for $\displaystyle b(p)\sum\limits_{\alpha\in \Lambda_p^{\lfloor\frac{p}{2}\rfloor}}\hat a^{2\alpha} \hat b^{\bar \alpha}$ in both cases, $p$ even or odd, we obtain that 
$$b(p)\sum\limits_{j=2}^{\lfloor\frac{p}{2}\rfloor}(-1)^j\sum\limits_{\alpha\in \Lambda_p^j}\hat a^{2\alpha} \hat b^{\bar \alpha} +R_p(a;b)=
\sum\limits_{j=2}^{\lfloor\frac{p+1}{2}\rfloor}(-1)^j \sum\limits_{\alpha \in \Lambda^{j}_{p+1}}  a^{2\alpha}b^
{\bar \alpha}$$ 
and hence that $Q_p(\hat a;\hat b)=Q_{p+1}(a;b)$.

If we consider now $a,c\in \ell(\KK^*;p,r)$, $b\in \ell(\KK;p,r)$, then  $\gamma=\sqrt{\dfrac{c^{\pi_p}}{ra^{\pi_p}}}$ and hence, 
$$\begin{array}{rl}
q_{p,r}(a;b;c)=Q_p(a_\phi;b_\phi)=&\hspace{-.25cm}\displaystyle \dfrac{1}{2a_\phi^{\pi_p}}
\sum\limits_{j=0}^{\lfloor\frac{p}{2}\rfloor}(-1)^j\sum\limits_{\alpha\in \Lambda_p^j}a_\phi^{2\alpha} b_\phi^{\bar \alpha}\\[3ex]
=&\hspace{-.25cm}\displaystyle \dfrac{1}{2 a^{\pi_p}\phi(a,c)^{\pi_p}}
\sum\limits_{j=0}^{\lfloor\frac{p}{2}\rfloor}(-1)^j\sum\limits_{\alpha\in \Lambda_p^j}\gamma^{2\alpha_{p-1}-1}a^{2\alpha} b^{\bar \alpha}\phi(a,c)^{2\alpha+\bar \alpha}\\[3ex]
=&\hspace{-.25cm}\displaystyle \dfrac{1}{2 \phi(a,c)^{\pi_p}}\sqrt{\dfrac{r}{a^{\pi_p} c^{\pi_p}}}
\sum\limits_{j=0}^{\lfloor\frac{p}{2}\rfloor}(-1)^j\sum\limits_{\alpha\in \Lambda_p^j}\left(\dfrac{c^{\pi_p}}{ra^{\pi_p}}\right)^{\alpha_{p-1}}a^{2\alpha} b^{\bar \alpha}\phi(a,c)^{2\alpha+\bar \alpha}.\end{array}$$
The result follows taking into account that 
 for any $\alpha\in \Lambda_p^m$, $m=0,\ldots,\lfloor\frac{p}{2}\rfloor$ we get 
$$\left(\dfrac{c^{\pi_p}}{a^{\pi_p}}\right)^{\alpha_{p-1}}\phi(a,c)^{2\alpha+\bar \alpha}=\phi(a;c)^{\pi_p}\dfrac{c^\alpha}{a^\alpha}.$$
\qed 
\vspace{.5cc}

Taking into account that from Lemma \ref{qp:cha}, $\ell(\CC;p,r)\subset \ell(\CC;np,r^n)$, for any $p\in \NN^*$  any $r\in \KK^*$ and any $n\in \NN^*$, then $q_{np,r^n}(a;b;c)$ has sense when $a,c\in \ell(\CC^*;p,r)$ and  $b\in \ell(\CC;p,r)$. In fact, we have the following relation between the Floquet functions $q_{np,r^n}$ and $q_{p,r}$

 \begin{propo}
Given $p\in \NN^*$ and  $r\in \CC^*$ then for any $a,c\in \ell(\CC^*;p,r)$ and  $b\in \ell(\CC;p,r)$, then for any $n\in \NN^*$ 
we have that 
$$q_{np,r^n}(a;b;c)=T_n\big(q_{p,r}(a;b;c)\big).$$
\end{propo}
\proof First, suppose that $a,c\in \CC^*$ and $b\in \CC$; that is, that the equation \eqref{equation} has constant coefficients. Then, $p=r=1$ and from Theorem \ref{Ff} and Proposition \ref{cardinal}, for any $n\in \NN^*$ we have 
$$q_{n,1}(a;b;c)=\dfrac{1}{2 }\sqrt{\dfrac{1}{a^nc^n}}
\sum\limits_{j=0}^{\lfloor\frac{n}{2}\rfloor}(-1)^ja^jb^{n-2j}c^j|\Lambda_n^j|=\dfrac{1}{2 } 
\sum\limits_{j=0}^{\lfloor\frac{n}{2}\rfloor}(-1)^j\dfrac{n}{n-j}{n-j\choose j}\Big(\dfrac{b}{\sqrt{ac}}\Big)^{n-2j}=T_n(q),$$
where $q=\dfrac{b}{2\sqrt{ac}}=q_{1,1}(a;b;c)$.

Assume now that  $a,c\in \ell(\CC^*;p,r)$ and  $b\in \ell(\CC;p,r)$ and consider $z\in \ell(\KK)$ a solution of the Equation 
$$a(k)z(k+1)-b(k)z(k)+c(k-1)z(k-1)=0,\hspace{.25cm}k\in \ZZ.$$

If  $\gamma=\big(r\phi(a,c)(p)\big)^{-\frac{1}{2}}$, then from Theorem \ref{Ch:sys}, for any $m\in \ZZ$, the sequence $v(k)=\gamma^{-k}z_{p,m}(k)$ is a solution of the Chebyshev equation with parameter $q_{p,r}(a;b;c)$.

On the other hand, since $\phi(a,c)(np)=\phi(a,c)(p)^n$ we have that $\big(r^n\phi(a,c)(np)\big)^{-\frac{1}{2}}=\gamma^n$ and hence, Theorem  \ref{Ch:sys} also implies that for any $m\in \ZZ$, the sequence $w(k)=\gamma^{-nk}z_{np,m}(k)$ is a solution of the Chebyshev equation with parameter $q_{np,r^n}(a;b;c)$.

Finally, since $z_{np,m}(k)=z(knp+m)=z_{p,m}(nk)$, we have that $w(k)=v(nk)=v_{n,0}(k)$. Therefore, the result follows by applying the first part of this proof. \qed
\vspace{.5cc}

The particularization of the above Proposition to Chebyshev equations, leads to the following generalization of Lemma \ref{Ch:odd-even}.

\begin{corollary}
\label{Fc} 
Given  $z\in \ell(\CC)$ a Chebyshev sequence with parameter $q\in \CC$, 
then for any $n,m\in \NN^*$, the subsequence  $z_{n,m}$ is a  Chebyshev sequence  with parameter $T_n(q)$. \end{corollary} 

Applying the above result together with  the identity \eqref{comb}, for any $n,m\in \NN^*$ and any $k\in \ZZ$ we have the following relations between Chebyshev polynomials of first and second kind:
$$\begin{array}{rl}
T_{kn+m}(x)=&\hspace{-.25cm}T_m(x)U_k\big(T_n(x)\big)-T_{m-n}(x)U_{k-1}\big(T_n(x)\big),\\[1ex]
U_{kn+m}(x)=&\hspace{-.25cm}U_m(x)U_k\big(T_n(x)\big)-U_{m-n}(x)U_{k-1}\big(T_n(x)\big).
\end{array}$$

Taking $k=1$  at both identities  we obtain the well known relations, see \cite{MH03}
\begin{equation}
\label{com:first} 
2T_m(x)T_n(x)=T_{n+m}(x)+T_{m-n}(x)\hspace{.25cm}\hbox{and}\hspace{.25cm}2U_m(x)T_n(x)=U_{n+m}(x)+U_{m-n}(x),
\end{equation}

Taking  $m=0$  in the first equation and  $n=2m+2$ in the second one,  we obtain the following generalizations of  the first identity in both \eqref{doubling:second} and   \eqref{doubling:first} 
\begin{equation}
\label{com:first} 
T_{kn}(x)=T_k\big(T_n(x)\big)\hspace{.25cm}\hbox{and}\hspace{.25cm}U_{2k(m+1)+m}(x)=U_{m}(x)W_k\big(T_{2(m+1)}(x)\big).
\end{equation}
Finally, taking $n=2m$ in the first equation and $m=n-1$ in the second one, we obtain the following generalizations of the second identity in \eqref{doubling:first} and in \eqref{doubling:second}, respectively
\begin{equation}
\label{doubling:general}
T_{m(2k+1)}(x)= T_m(x)V_k\big(T_{2m}(x)\big)\hspace{.25cm}\hbox{and}\hspace{.25cm}U_{(k+1)n-1}(x)= U_{n-1}(x)U_k\big(T_n(x)\big).
\end{equation}

We end this paper,  analyzing when  a second order difference equation with quasi-periodic coefficients with period $p$ has also quasi-periodic solutions with the same period.  So, given $p\in \NN^*$, $r\in \KK^*$ and the sequences  $a,c\in \ell(\KK^*;p,r)$, $b\in \ell(\KK;p,r)$, if $z\in \ell(\KK;p,\hat r)$ is a solution of Equation (\ref{equation}), from Lemma \ref{qp:cha} we know that $z_{p,m}\in \ell(\KK;1,\hat r)$. Moreover,  Theorem \ref{Ch:sys} establishes that $z_{p,m}(k)=\gamma^k v(k)$, where $\gamma=\sqrt{\dfrac{c^{\pi_p}}{ra^{\pi_p}}}$ and $v(k)$ is a solution of the Chebyshev equation with parameter $q_{p,r}(a;b;c)$. Therefore, $z\in \ell(\KK;p,\hat r)$ iff $v\in \ell(\KK;1,\hat r\gamma^{-1})$.  So, applying Proposition \ref{Cheb:periodic} we have the following characterization.

\begin{theorem} 
\label{qperiodic:solutions}
Given $p\in \NN^*$ and  $r\in \KK^*$ then for any $a,c\in \ell(\KK^*;p,r)$ and  $b\in \ell(\KK;p,r)$, the difference equation 
$$a(k)z(k+1)-b(k)z(k)+c(k-1)z(k-1)=0,\hspace{.25cm}k\in \ZZ,$$
has quasi-periodic solutions with period $p$ and ratio $\hat  r\in \KK^*$ iff 
$$\hat r a^{\pi_p}+(r\hat r)^{-1}c^{\pi_p}=
\sum\limits_{j=0}^{\lfloor\frac{p}{2}\rfloor}(-1)^j\sum\limits_{\alpha\in \Lambda_p^j}r^{-\alpha_{p-1}}a^{\alpha} b^{\bar \alpha}c ^{ \alpha}.$$
\end{theorem}
 
When $r=\hat r=1$, the above result establishes the necessary and sufficient condition for the existence of periodic solutions for difference equations with periodic coefficients, which represent a fully generalization of Corollary \ref{floquet:constant}.

\begin{corollary}[Floquet] If  $a,c\in \ell(\KK^*;p)$ and  $b\in \ell(\KK;p)$, the difference equation 
$$a(k)z(k+1)-b(k)z(k)+c(k-1)z(k-1)=0,\hspace{.25cm}k\in \ZZ,$$
has periodic solutions with period $p$ iff 
$$
a^{\pi_p}+c^{\pi_p}=
\sum\limits_{j=0}^{\lfloor\frac{p}{2}\rfloor}(-1)^j\sum\limits_{\alpha\in \Lambda_p^j}a^{\alpha} b^{\bar \alpha}c ^{ \alpha}.
$$
\end{corollary}
The reader could compare the condition given in the above corollary  with the same result in \cite[Corollary 2.9.2]{A00}. 
 \vspace{.5cc}

 {\bf Acknowledgments.}
This work has been partly supported by the Spanish Research Council (Comisi\'on Interministerial de Ciencia y Tecnolog\'\i a,)
under project MTM2011-28800-C02-02.


\vspace{.5cm}

\begin{thebibliography}{99}

\bibitem{A00}
Agarwal, R.P., {\it Difference equations and inequalities},
Marcel Dekker, 2000.

\bibitem{ABD05} Aharonov, D., Beardon A., Driver, K.
\newblock Fibonacci, Chebyshev, and Orthogonal Polynomials.
\newblock  {\em Amer. Math. Monthly}, {\bf 112} (2005), 612-630.


%
%
%
%
%
%
%
%
%


%

%

\bibitem{M97} Mallik, R.K.
\newblock On the Solution of a Second Order Linear
Homogeneous Difference Equation with Variable
Coefficients.
\newblock  \emph{J. Math. Anal. Appl.}, {\bf 215} (1997), 32-47.


\bibitem{MH03} Mason, J.C., Handscomb, D.C.
\newblock\emph{Chebyshev Polynomials}.
\newblock Chapman \& Hall/CRC, 2003.


%
%
%
%
%
%
%
%
%
%


\end{thebibliography}
\end{document}